\documentclass[mathematics,article,accept,moreauthors,pdftex]{mdpi} 

\usepackage{doi}

\DeclareMathOperator{\dv}{div}
\DeclareMathOperator*{\indmax}{indmax}

\newcommand{\eps}{\varepsilon}
\newcommand{\way}{{\{way\}}}
\newcommand{\nmax}[1]{n^{\{#1\}}_{\max}}
\newcommand{\lmax}[1]{l^{\{#1\}}_{\max}}
\newcommand{\dtmax}[1]{\Delta t^{\{#1\}}_{\max}}
\newcommand{\dmax}[1]{d^{\{#1\}}_{\max}}

\newcommand{\reff}[1]{(\ref{#1})}
\newcommand{\tr}{^{\scriptsize\mbox{T}}}
\newcommand{\disp}{\displaystyle}
\newcommand{\prodser}[3]{\disp\sum_{j=0}^{#1}\langle #3Q_{#2}#3\lam_j,#3\lam_{#1-j}\rangle}
\newcommand{\nl}{\medskip\\}
\newcommand{\lam}{\Upsilon}
\newcommand{\nr}[1]{\left| #1 \right|}

\firstpage{1} 
\pubvolume{9}
\issuenum{17}
\articlenumber{2057}
\pubyear{2021}
\copyrightyear{2021}
\datereceived{18 July 2021} 
\dateaccepted{17 August 2021} 
\datepublished{26 August 2021} 
\hreflink{https://doi.org/10.3390/math9172057}

\Title{On the Poisson Stability to Study a Fourth-Order Dynamical System with Quadratic Nonlinearities}

\TitleCitation{On the Poisson Stability to Study a Fourth-Order Dynamical System with Quadratic Nonlinearities}

\Author{Alexander N. Pchelintsev 
	\orcidA{0000-0003-4136-1227}}

\AuthorNames{Alexander N. Pchelintsev}

\AuthorCitation{Pchelintsev, A.N.}

\address[1]{
	Tambov State Technical University, Department of Higher Mathematics,
	ul. Sovetskaya 106, 392000  Tambov,  Russia; pchelintsev.an@yandex.ru}

\abstract{This article discusses the search procedure for the Poincar\'e recurrences to classify solutions on an attractor of a fourth-order nonlinear dynamical system using a previously developed high-precision numerical method. For the resulting limiting solution, the Lyapunov exponents are calculated using the modified Benettin's algorithm to study the stability of the found regime and confirm the type of attractor.}



\begin{document}
\end{paracol}

\section{Introduction}

As it is known \cite{SprottBook}, the calculation of the Lyapunov exponents is used for the classification
of attractors of dynamical systems. The combinations of signs of such values determine the attractor type:
an equilibrium position, limit cycle, multidimensional torus or strange attractor.
The dynamics of the corresponding solutions are called static, periodic, quasiperiodic, or chaotic.
The computational procedure for determining the Lyapunov exponents is mainly based on the Benettin's
algorithm \cite{Tancredi}. However, in \cite{LoziPchelintsev, LozPogPchel} the points
of limiting solutions were investigated for the Poisson stability, which made it possible to understand
whether we have a quasiperiodic or chaotic regime.
In \cite{LozPogPchel}, it was done for the Jafari--Sprott system \cite{Jafari}.

It is noteworthy \cite{Nemytskii} that a point $y$ of the phase space is called positively Poisson stable
($P^+$), if for any neighborhood $U$ of the point $y$ and for any $T_P> 0$ there is a time
value $t\ge T_P$ for which the trajectory of the dynamical system enters into the neighborhood $U$.
Similarly, if there is $t\le -T_P$ such that the trajectory enters into the neighborhood $U$, then
the point $y$ is negatively Poisson stable ($P^-$). A point is called Poisson stable,
if it is $P^+$ and $P^-$ -- stable.

The boundedness of the limiting solutions of dissipative systems implies
\cite{Nemytskii,Dzyuba1,Dzyuba2,Dzyuba3,Kuznetsov} that any steady-state oscillation mode is described by Poisson-stable trajectories. This also applies to dynamical chaos. A trajectory different from the equilibrium position is said to be Poisson stable if it verifies the property of
returning in an arbitrarily small $\eps$-neighborhood of each of its points an infinite number of times.
Such returns are called the Poincar\'e recurrences. In \cite[p. 146]{Kuznetsov}, the author notes that "the
study of the statistics of the Poincar\'e recurrences is a powerful tool for the analysis and classification
of dynamic modes. Apparently, the potentialities of this approach have not yet been fully exhausted
in modern nonlinear dynamics". For example, the returns
follow one another regularly for quasiperiodic regimes.
Then \cite[p. 145]{Kuznetsov} "the dynamic chaos is a situation when the Poincar\'e recurrences to the
$\eps$-neighborhood of the initial point do not show regularity, the time interval between two
successive returns turns out to be different each time, and some statistical distribution of the
times of return is arised" \cite{AfraimovichLin,AfraimovichZaslavsky}.
An example of the Poincar\'e recurrences analysis based on the Kac's theorem
\cite{Kac,Penne,Afraimovich,Anishchenko1,Anishchenko2,Anishchenko3}, for a discrete dynamical
system with a chaotic non-hyperbolic attractor is given in \cite{Anishchenko1}.

As it is known, for the most of dissipative systems, the formulas for a general solution of
the system in a class of any known functions with well-studied properties have not yet been found.
Therefore, numerical methods are used. The use of classical numerical methods (such as
Euler method, Runge-Kutta 4th order method, Adams methods, etc.) for constructing approximate
solutions in attractors of dynamical systems leads to significant errors over large sections of
time due to the instability of the studied chaotic regimes.
In general, the problem of numerical modeling with control of the accuracy of obtained solution
and the choice of a platform for computer implementation is relevant today \cite{Nepomuceno},
since small errors introduced at each integration step cause exponential divergence of close
trajectories. More recently, the numerical FGBFI method (Firmly Grounded Backward-Forward
Integration) has been proposed in \cite{LoziPchelintsev,LozPogPchel,Pchel2020}, based
on the power series method for dynamical systems with quadratic nonlinearities.
The main advantages of this method are as follows:
\begin{enumerate}
\item The recurrence relations are obtained for calculating the coefficients of the expansion
of solutions in a power series for any dynamical system with quadratic nonlinearities
\textit{in a general form};

\item The convergence of the power series is studied. A simple formula for calculations is
derived (in comparison with that \cite{Babadzanjanz} obtained in the known
literature) for calculating the length of the integration step \textit{in a general form};

\item The criteria for checking the accuracy of the approximate solution are obtained.
The control of the accuracy and configuration of the approximate solution of a dynamical
system uses forward and backward time, which makes the numerical method reliable (degrees
of piecewise polynomials, the value of the maximum integration step, etc.);

\item The FGBFI method allows to construct high-precision approximations to non-extendable
solutions of a system of autonomous differential equations with a quadratic right-hand side.
Like for instance the system of the form
$$
  \dot{x}=1+x^2.
$$
In this case, the numerical solution computed with FGBFI will never cross the asymptote and will approach it
arbitrarily close.
\end{enumerate}

Thus, it is a high-precision method for constructing the trajectory arc of the system for
any time interval. This makes it possible to track the Poincar\'e recurrences to any neighborhood of the
trajectory points, since the resulting computational error can be smaller than the
radius of the monitored neighborhoods.

This article considers a fourth-order dissipative system \cite{Dong2019}, in which,
\textit{in the opinion of its authors}, there is a hyperchaos. The goals of the research are:
(1) to find the Poincar\'e recurrences in the attractor of this system for the approximate solutions
obtained by the FGBFI method and to analyze the statistics of return times; (2) to apply the modified
Benettin's algorithm \cite{Pchel2020} to refine the computation of the Lyapunov exponents.

\section{Bohr's almost periodic functions}

It is noteworthy \cite[pp. 368, 418, 419]{Demidovich} that the function (trajectory of a dynamical system)
${F(t)\in\mathbb R^n}$
is called Bohr's almost periodic in the sense that if for any ${\eps>0}$ there exists a relatively
dense set of almost periods $\tau_F=\tau_F(\eps)$ of the function $F(t)$ with $\eps$-accuracy,
i.e. there is the positive number $L=L(\eps)$ such that any interval $[\alpha;\alpha+L]$ contains
at least one number $\tau_F$ for which the inequality
$$
  \nr{F(t+\tau_F)-F(t)}<\eps
$$
holds at $t\in\mathbb R$. Note that for an almost periodic function (and, in general, for a Poisson
stable trajectory), different from a periodic function, with decreasing number $\eps$
the number $L(\eps)$ must increase indefinitely, otherwise the function would be periodic.

Usually (see, for example, \cite{Kuznetsov}) in nonlinear dynamics, almost periodic functions
are called quasiperiodic.

Note an important theorem in the theory of dynamical systems \cite[p. 147]{Kuznetsov}:
if a non-periodic trajectory is Poisson and Lyapunov stable, then it is almost periodic.

Since we will be investigating the Poincar\'e recurrences, it is noteworthy that for an almost
periodic trajectory (as follows from the above definition), such returns must form a relatively
dense set. An example of a sequence that is not a relatively dense set is:
$$
  0,\:1^2,\:2^2,\:3^2,\:\ldots,
$$
since
$$
  \sup_q |(q+1)^2-q^2|=\infty,
$$
i.e. in the given example, the distance between neighboring elements grows with the growth of $q$.
A similar situation was observed in a numerical experiment for the Poincar\'e recurrences in the case of
a chaotic solution of the Chen system \cite{LoziPchelintsev}.

\section{Dissipative system of the fourth order}

Let us consider a fourth-order system \cite{Dong2019}
\equation
   \left\{\begin{array}{l}
      \dot{x}_1=ax_2-ax_1-ex_4,\\
      \dot{x}_2=bx_1-x_2-x_1x_3-fx_4,\\
      \dot{x}_3=x_1x_2-cx_3,\\
      \dot{x}_4=kx_2x_3-dx_4,
   \end{array}\right.\label{dynsys}
\endequation
where $a$, $b$, $c$, $d$, $e$, $f$ and $k$ are positive parameters. The divergence of the
vector field $G$ defined by the right-hand side of the system \reff{dynsys} is equal to
$$
  \dv G=-(a+1+c+d)<0.
$$
Then the system \reff{dynsys} is dissipative. Hence, there is the ball $B_a\subset\mathbb R^4$,
into which each trajectory of the system \reff{dynsys} is immersed forever, after a while.
Therefore, there is a limit set (attractor) such that all trajectories of the
dynamical system are attracted when $t\rightarrow\infty$ \cite{Nemytskii}.

Let us investigate the dynamics of the system \reff{dynsys} for the values of the parameters
$a=7$, $b=50$, $c=3$, $d=10$, $e=5$, $f=5$ and $k=1.5$. In this case, it is possible to determine
the position and radius of the ball $B_a$ from the data obtained in \cite{Dong2019}.
Also, in this article the Lyapunov exponents were determined: $LE_1=2.1040$, $LE_2=0.3563$,
$LE_3=0$ and $LE_4=-23.4508$, indicating a hyperchaos in the system \reff{dynsys}.
Therefore, its solutions are highly unstable, which requires the use of special numerical
procedures for large time intervals.

Let us consider the application of the FGBFI method for constructing approximations to the
solutions of the system \reff{dynsys} using high-precision calculations.

\section{The FGBFI method}
\label{algo_FGBFI}

Following \cite{LozPogPchel}, we will rewrite the system \reff{dynsys} as
\equation
   \dot{X}=AX+\Phi(X),\label{polsys}
\endequation
where
$$
\begin{array}{c}
   X(t)=\left[x_1(t)\:\:x_2(t)\:\:x_3(t)\:\:x_4(t)\right]\tr,\:
   \Phi(X)=\left[\varphi_1(X)\:\:\varphi_2(X)\:\:\varphi_3(X)\:\:\varphi_4(X)\right]\tr,\nl
   \varphi_p(X)=\langle Q_p X,X \rangle,\:\:\:p=\overline{1;4},
   \:\:\:A=\left[
           \begin{array}{rrrr}
              -a &  a &  0 & -e\\
               b & -1 &  0 & -f\\
               0 &  0 & -c &  0\\
               0 &  0 &  0 & -d
           \end{array}
           \right],\nl
         Q_1=\textbf{0},\:\:\:
         Q_2=\left[
             \begin{array}{rrrr}
                0 & 0 & -1 & 0\\
                0 & 0 &  0 & 0\\
                0 & 0 &  0 & 0\\
                0 & 0 &  0 & 0
             \end{array}
             \right],
\end{array}
$$
$$
\begin{array}{c}
         Q_3=\left[
             \begin{array}{rrrr}
                0 & 1 & 0 & 0\\
                0 & 0 & 0 & 0\\
                0 & 0 & 0 & 0\\
                0 & 0 & 0 & 0
             \end{array}
             \right],\:\:\:
         Q_4=\left[
             \begin{array}{rrrr}
                0 & 0 & 0 & 0\\
                0 & 0 & k & 0\\
                0 & 0 & 0 & 0\\
                0 & 0 & 0 & 0
             \end{array}
             \right].
\end{array}
$$

Let us expand the solution
\equation
   X(t)=\sum_{i=0}^\infty \lam_i t^i,\label{powser}
\endequation
where $X(0)=\lam_0$ is an initial condition for the system \reff{polsys}, the vector $\lam_0$ is given,
$$
  \lam_i=\left[\upsilon_{i,1}\:\:\upsilon_{i,2}\:\:\upsilon_{i,3}\:\:\upsilon_{i,4}\right]\tr,\:\:
  \upsilon_{i,p}\in\mathbb R,\:\:p=\overline{1;4}.
$$

Let
$$
  \Psi_i=\left[\psi_{i,1}\:\:\psi_{i,2}\:\:\psi_{i,3}\:\:\psi_{i,4}\right]\tr
$$
and
$$
  \psi_{i,p}=\prodser{i}{p}{},\:\:p=\overline{1;4}.
$$

The recurrence relations for calculating the coefficients of the power series have
the form
\cite{LozPogPchel}
\equation
   \lam_j=\dfrac{A\lam_{j-1}+\Psi_{j-1}}{j}\label{cser}
\endequation
for $j\in\mathbb N$.

Since the criteria for checking the accuracy of the obtained numerical solution require
repeated calculations in forward and backward time, we need to have a guaranteed estimate
of the convergence interval of the power series \reff{powser} for a given vector $\lam_0$.
In \cite{LozPogPchel}, a theorem is proved on the estimate of this interval for
systems with a quadratic right-hand side.

The norms are calculated:
$$
  \begin{array}{c}
     \|A\|=\|A\|_1=\max\{a+b,\:a+1,\:c,\:e+f+d\}=57,\nl
     \|Q_1\|=0,\:\:\:
     \|Q_2\|=\|Q_3\|=1,\:\:\:
     \|Q_4\|=k,\nl
     \disp
     \mu=4\cdot\max_{p=\overline{1;4}}\|Q_p\|=4\cdot1.5=6.
  \end{array}
$$
Next, two auxiliary numbers are calculated as functions of $\lam_0$:
$$
  \begin{array}{c}
     h_1(\lam_0)=\|\lam_0\|=\disp\sum_{p=1}^4|\upsilon_{0,p}|,\nl
     h_2(\lam_0)=\left\{\begin{array}{l}
                           \mu h_1^2(\lam_0)+(\|A\|+2\mu)h_1(\lam_0),
                           \:\:\mbox{if}\:\:h_1>1,\\
                           \|A\|+\mu\:\:\mbox{otherwise},
                        \end{array}
                 \right.
  \end{array}
$$
or
$$
  h_2(\lam_0)=\left\{\begin{array}{l}
                        6\cdot\|\lam_0\|^2+69\cdot\|\lam_0\|,
                          \:\:\mbox{if}\:\:\|\lam_0\|>1,\\
                        63\:\:\mbox{otherwise}.
                     \end{array}
              \right.
$$
Then the series \reff{powser} converges for $t\in[-\tau,\tau]$, where
\equation
   \tau=\tau(\lam_0)=\dfrac{1}{h_2(\lam_0)+\delta},\label{tauser}
\endequation
$\delta$ is any positive number (can be taken very small).

It can be seen from the formula \reff{tauser} that the length of the
convergence interval depends on the choice of the initial condition for the
system \reff{polsys}. This fact is confirmed by computational experiments,
for example, for the Lorenz system \cite{LorPchel}.

It is worth noting that a similar approach can be applied to systems with
cubic nonlinearities, but when the derivative on the left side of the system is
of the second (or higher) order. For example, in \cite{PchelAhmad}
this was done for the Duffing equation.

Usually, in calculations, many researchers work with subnormal real numbers of
single or double precision, presented in the IEEE 754 format \cite{Overton}.
The main drawback here is the fixed accuracy of the representation of real numbers,
which may not allow us numerically constructing approximations to unstable solutions
of systems of differential equations on the large time intervals.
Therefore, for high-precision calculations, the GNU MPFR library \cite{MPFR,Fousse} is used.

Next, consider the algorithm \cite{Pchel2020} for the numerical construction of
an approximate solution of a dynamical system in forward and backward time:

\begin{enumerate}
   \item \textbf{Set} the number $b_m$ of bits under the mantissa of a real
   number and the precision $\eps_{pw}$ for estimate of the common term of
   power series. Note that $b_m$ defines the machine epsilon $\eps_m$.
   Let us choose $b_m$ so that the precision of representation of the real
   number is with a margin, i.e.
   $$
     \eps_m\ll\eps_{pw};
   $$
   \item $t:=0$;
   \item \textbf{Set} $X(0)\in B_a$ for the system \reff{polsys} and value the number $way$
   that determines the direction in time: $way=1$ is gone forward in time, $way=-1$ is
   gone backward in time;
   \item \textbf{Set} $T$: the length of the time interval on which the numerical integration
   will be performed;
   \item $ended:=\mbox{{\ttfamily false}}$;
   \item\label{stepgoto} \textbf{Calculate} the integration step $\Delta t:=\tau(X(0))$ by the formula
        \reff{tauser};
   \item \textbf{If} $\Delta t>T-t$, \textbf{then} $\Delta t:=T-t$, $t:=T$
   
   \textbf{Else} $t:=t+\Delta t$;
   \item $\Delta t:=way\cdot\Delta t$;
   \item \textbf{Let} $\lam_0:=X(0)$;
   \item \textbf{Calculate} the approximate value $X(\Delta t)$ by summing the terms of the series
         \reff{powser} to such a value $i$, where the following inequality holds
         \equation
            \left\|\lam_i\right\|\cdot|\Delta t|^i<\eps_{pw};\label{epssr}
         \endequation
   \item \textbf{Print} $way\cdot t,\:X(\Delta t)$;
   \item \textbf{If} $X(\Delta t)\notin B_a$, \textbf{then} we got out the compact $B_a$. Then
         \textbf{write} {\sffamily "Decrease the value $\eps_{pw}$ and/or $\eps_m$"};
         $ended:=\mbox{{\ttfamily true}}$;
   \item \textbf{If} $t=T$, \textbf{then} $ended:=\mbox{{\ttfamily true}}$;
   \item \textbf{If} $ended$, \textbf{then} \textbf{terminate} the algorithm;
   
   \item $X(0):=X(\Delta t)$;
   \item \textbf{Goto} step \ref{stepgoto}.
\end{enumerate}

As can be seen, the algorithm is universal in the direction of numerical integration in time.

Next, we consider the criteria for verifying the accuracy of the resulting solution.
For this we assume
$$
  \begin{array}{c}
     \Omega=[t_0,t_1]\cup[t_1,t_2]\cup\ldots\cup[t_{N-1},t_N],\nl
     t_0=0,\:\:t_N=way\cdot T,\nl
     \Delta t_l^\way=t_l-t_{l-1}=\tau\left(X_{l-1}\left(\Delta t_{l-1}^\way\right)\right),
     \:\:X_0\equiv X(0),\:\:l=\overline{1;N^\way},
  \end{array}
$$
where $l$ is an index of the interval $[t_{l-1},t_l]$, on which the series \reff{powser} is converged,
$N^\way$ is a number of such intervals for the direction $way$ on time. Since on each interval $[t_{l-1},t_l]$
during calculations we truncate the series \reff{powser} to some polynomial $X_l(t)$ approximating the
solution of the system \reff{polsys} on it (i.e. locally), then we denote by $n_l^{\{way\}}$ the degree
of the polynomial $X_l(t)$. Then we introduce the notation:
\equation
   \begin{array}{c}
      \disp\nmax{way}=\max_l n^\way_l,\:\:
      \lmax{way}=\indmax_l n^\way_l,\:\:
      \disp\dtmax{way}=way\cdot\max_l \left|\Delta t^\way_l\right|,\nl
      \dmax{way}=\indmax_l \left|\Delta t^\way_l\right|.
   \end{array}
\label{nums}
\endequation

In \cite{LozPogPchel}, the following criteria for checking the accuracy of
the obtained solution were proposed:
\begin{enumerate}
   \item The accuracy $\eps_a$ of approximation at $way=1$ is a frequently used criterion
         in applications of numerical methods for solving differential equations. When the
         inequality \reff{epssr} is true, it is necessary to increase the degrees of all
         polynomials $n^{\{1\}}_l$, obtaining the next approximation, and compare the
         distance $\delta_a$ between the obtained approximate solutions on the interval $[0,T]$
         with the value $\eps_a$. If $\delta_a>\eps_a$, then we increase the powers of $n_l$,
         otherwise we have to use the obtained solution;
   \item The radius $\eps_R$ of the neighborhood of the initial point, to which the approximate
         solution should return in backward time is another criterion. In other words, we need
         to select the accuracy of $\eps_{pw}$ so that the following inequality holds
         $$
           \nr{X_{\widehat{N}}^{\{-1\}}\left(\Delta t_{\widehat{N}}^{\{-1\}}\right)-X(0)}<\eps_R,
         $$
         where $\widehat{N}=N^{\{-1\}}$. Note here that we do not know the exact solution to the system
         \reff{polsys}, but we do know its initial point. Then we can set how many digits are in each
         component of the vector $X_{\widehat{N}}^{\{-1\}}\left(\Delta t_{\widehat{N}}^{\{-1\}}\right)$
         after the point must match the digits of the corresponding components $X(0)$. The problem is
         that a large amount of computation is required (due to repetitions of the forward and backward
         traverse in time). The solutions of the system \reff{polsys}, as a rule, are highly unstable in
         the backward time: they immediately leave from the attractor, because in our calculations, we are
         near it, no matter how accurate $\eps_{pw}$ is taken. Therefore, in the above algorithm, we control
         the finding of the trajectory within the boundaries of the set $B_a$;
   \item The comparison of configurations of approximate solutions in forward and backward time is
         necessary to determine the numbers \reff{nums}, describing approximate solutions. Next, check
         $$
           \begin{array}{c}
              N\approx\widehat{N},\:\nmax{1}\approx\nmax{-1},\:
              t_{\lmax{1}}+\left|t_{\lmax{-1}}\right|\approx T,\nl
              \dmax{1}+\dmax{-1}\approx N,\:
              t_{\dmax{1}}+\left|t_{\dmax{-1}}\right|\approx T.
            \end{array}
         $$
         Unlike criteria 1 and 2, we control here the arguments of the approximate solutions.
\end{enumerate}

\section{Analysis of the Poincar\'e recurrences for a fourth-order system}

\begin{table}[!]
	\footnotesize
	\caption{The Poincar\'e recurrences for the point \reff{init_cond2}, $\Delta t_P=10^{-4}$, $T=15$ for reaching
		the attractor.}
	\centering\vspace{10pt}
	{\begin{tabular}{|c|c|}\hline
			Time moment $t_{k^*}$ & Value $d_{k^*}$\\\hline
			0.3655 & 0.0423764\\\hline
			0.7310 & 0.0157840\\\hline
			1.0965 & 0.0323213\\\hline
			1.4620 & 0.0196461\\\hline
			1.8275 & 0.0173865\\\hline
			2.1930 & 0.0209083\\\hline
			2.5585 & 0.0135900\\\hline
			2.9240 & 0.0206764\\\hline
			3.2895 & 0.0194286\\\hline
			3.6550 & 0.0197196\\\hline
			4.0205 & 0.0237665\\\hline
			4.3860 & 0.0213908\\\hline
			4.7515 & 0.0263917\\\hline
			5.1170 & 0.0259931\\\hline
			5.4825 & 0.0288078\\\hline
			5.8480 & 0.0311426\\\hline
			6.2135 & 0.0320319\\\hline
			6.5790 & 0.0356643\\\hline
			6.9445 & 0.0362604\\\hline
			7.3101 & 0.0340862\\\hline
			7.6756 & 0.0325096\\\hline
			8.0411 & 0.0305410\\\hline
			8.4066 & 0.0280697\\\hline
			8.7721 & 0.0268155\\\hline
			9.1376 & 0.0243143\\\hline
			9.5031 & 0.0231314\\\hline
			9.8686 & 0.0211890\\\hline
	\end{tabular}}
	\label{tab:1}
\end{table}

\begin{table}[!]
	\footnotesize
	\caption{The Poincar\'e recurrences for the point \reff{init_cond3}, $\Delta t_P=10^{-4}$, $T=40$ for reaching
		the attractor.}
	\centering\vspace{10pt}
	{\begin{tabular}{|c|c|}\hline
			Time moment $t_{k^*}$ & Value $d_{k^*}$\\\hline
			0.3655 & 0.00133861\\\hline
			0.7310 & 0.00267728\\\hline
			1.0965 & 0.00401586\\\hline
			1.4620 & 0.00535439\\\hline
			1.8275 & 0.00669292\\\hline
			2.1930 & 0.00803134\\\hline
			2.5585 & 0.00936979\\\hline
			2.9240 & 0.01070810\\\hline
			3.2895 & 0.01204650\\\hline
			3.6550 & 0.01338470\\\hline
			4.0205 & 0.01472300\\\hline
			4.3860 & 0.01606120\\\hline
			4.7515 & 0.01739930\\\hline
			5.1171 & 0.01774010\\\hline
			5.4826 & 0.01640090\\\hline
			5.8481 & 0.01506170\\\hline
			6.2136 & 0.01372250\\\hline
			6.5791 & 0.01238330\\\hline
			6.9446 & 0.01104430\\\hline
			7.3101 & 0.00970522\\\hline
			7.6756 & 0.00836622\\\hline
			8.0411 & 0.00702727\\\hline
			8.4066 & 0.00568835\\\hline
			8.7721 & 0.00434949\\\hline
			9.1376 & 0.00301067\\\hline
			9.5031 & 0.00167189\\\hline
			9.8686 & 0.00033316\\\hline
	\end{tabular}}
	\label{tab:2}
\end{table}

\begin{table}[!]
	\footnotesize
	\caption{The Poincar\'e recurrences for the point \reff{init_cond3}, $\Delta t_P=10^{-6}$, $T=40$ for reaching
		the attractor.}
	\centering\vspace{10pt}
	{\begin{tabular}{|c|c|}\hline
			Time moment $t_{k^*}$ & Value $d_{k^*} \cdot 10^{-5}$\\\hline
			0.365504 & 1.105410\\\hline
			0.731007 & 1.096510\\\hline
			1.096510 & 0.348186\\\hline
			1.462010 & 0.752419\\\hline
			1.827520 & 1.799740\\\hline
			2.193020 & 0.697982\\\hline
			2.558530 & 0.397935\\\hline
			2.924030 & 1.498770\\\hline
			3.289530 & 1.051780\\\hline
			3.655040 & 0.048264\\\hline
			4.020540 & 1.146250\\\hline
			4.386040 & 1.403530\\\hline
			4.751550 & 0.304245\\\hline
			5.117050 & 0.793828\\\hline
			5.482550 & 1.754720\\\hline
			5.848060 & 0.656571\\\hline
			6.213560 & 0.442487\\\hline
			6.579070 & 1.541010\\\hline
			6.944570 & 1.008280\\\hline
			7.310070 & 0.091034\\\hline
			7.675580 & 1.189060\\\hline
			8.041080 & 1.359960\\\hline
			8.406580 & 0.261577\\\hline
			8.772090 & 0.837369\\\hline
			9.137590 & 1.711900\\\hline
			9.503100 & 0.613282\\\hline
			9.868600 & 0.485547\\\hline
	\end{tabular}}
	\label{tab:3}
\end{table}

In \cite[p. 3221]{Dong2019}, the authors indicated that there is a hyperchaos
for the considered parameters of the system \reff{dynsys}. For example, the authors found it
for the initial condition
\equation
  x_1(0)=10,\:\:x_2(0)=50\cdot\sin 10\approx-27.2011,\:\:x_3(0)=10,\:\:x_4(0)=10.\label{init_cond1}
\endequation
Let us verify it.

Using the numerical method FGBFI, we construct the trajectory arc from the point \reff{init_cond1} on
the time interval $[0,T]$, $T=15$ getting close to the attractor. Let us indicate the parameters of
the numerical method:
$$
  b_m=128,\:\:\eps_m=5.87747\cdot10^{-39},\:\:\eps_{pw}=10^{-20}.
$$
We get the next point at the end time:
\equation
   X(15)=[
                   6.2355509634533960831\:\:\:
                   2.0140572482317481452\:\:\:
                  35.4929323328531102196\:\:
                 -43.5507482101916799734
         ]\tr.\label{init_cond2}
\endequation
\begin{paracol}{1}
Note that in the numerical experiment by to the first criterion for checking the accuracy,
decreasing the value of $\eps_{pw}$, the number of indicated correct signs in coordinates is
preserved.

Next, we take the point \reff{init_cond2} as the initial point and track the Poincar\'e recurrences for it.
In order to reduce the amount of calculations without moving backward in time (so as not
to take too small a value of $\eps_{pw}$), we look at the times when there is a return to the neighborhood
of the initial point for different values $\eps_{pw}$. This is necessary to control the accuracy of the
obtained points, choosing such the value $\eps_{pw}$ that these moments of time do not differ significantly.

To determine the Poincar\'e recurrences, we choose the small time step $\Delta t_P>0$. Let us divide the time
interval $[0,T_P]$, where $T_P$ is the end time at which the search for the Poincar\'e recurrences is terminated,
into intervals of length $\Delta t_P$ (let $N_P$ be the number of such intervals). Let us denote $X_k$
the coordinates of the $k$-th point of the considered trajectory arc, corresponding to the time moment
$$
  t_k=k\Delta t_P,\:\:k=\overline{1;N_P}.
$$

We introduce the distance
$$
  d_k=|X_k-X(0)|.
$$

To track the Poincar\'e recurrences, we fix such values $k=k^*$, moving in the trajectory arc, when there
are \textit{the local rapprochements with the point $X(0)$}, i.e.
$$
  d_{k^*-1}>d_{k^*}\:\:\:\&\:\:\:d_{k^*}<d_{k^*+1}\:\:\:\&\:\:\:d_{k^*}<1.
$$

For the point \reff{init_cond2} for $T_P=10$ and $\Delta t_P=10^{-4}$, the returns are shown in Tabl. \ref{tab:1}.
As can be seen from this table, the rapprochement with the initial point occurs at approximately
the same time (i.e., we have regularity), which suggests the idea of a limit cycle in the system \reff{dynsys}. Decrease in the value of $\Delta t_P$ does not cause significant changes in Tabl.
\ref{tab:1}.

In this case to show the absence of chaos in the system \reff{dynsys}, we take the time value $T=40$ (instead of
$T=15$, which was used to obtain the coordinates of the point \reff{init_cond2}) to reach the attractor. We get
the following point at the end time:
\equation
   X(40)=[
                   1.6321991613781496393\:\:\:
                   8.7300523565474285155\:\:\:
                  39.6961687172415982460\:\:\:
                  54.8461996449311966025
         ]\tr.\label{init_cond3}
\endequation
Similarly, for the point \reff{init_cond3}, the Poincar\'e recurrences shown in Tabl. \ref{tab:2} for $T_P=10$ and
$\Delta t_P=10^{-4}$. As we can be seen from this table, returns are getting closer. Therefore, we will decrease the
value $\Delta t_P$, setting it equal to $\Delta t_P=10^{-7}$. We get the returns shown in Tabl.
\ref{tab:3}.

\begin{figure}[!]
   \centering\includegraphics[width=\textwidth,height=9cm]{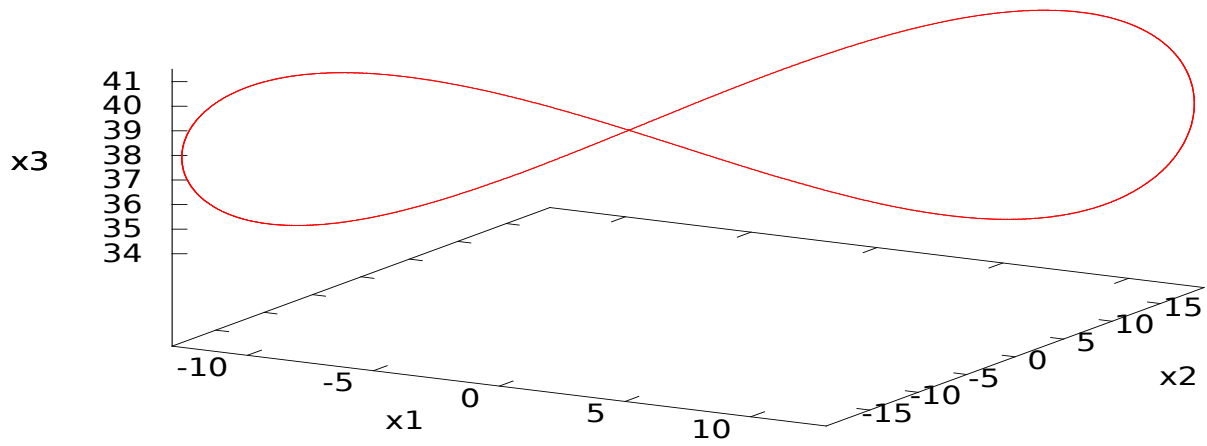}
   \caption{Projection of the cycle of the system \reff{dynsys} with the period $T_c \approx 0.3655$
   	        into the space ${x_1-x_2-x_3}$.}\label{fig:1}
\end{figure}

\begin{figure}[!]
   \centering\includegraphics[width=\textwidth,height=9cm]{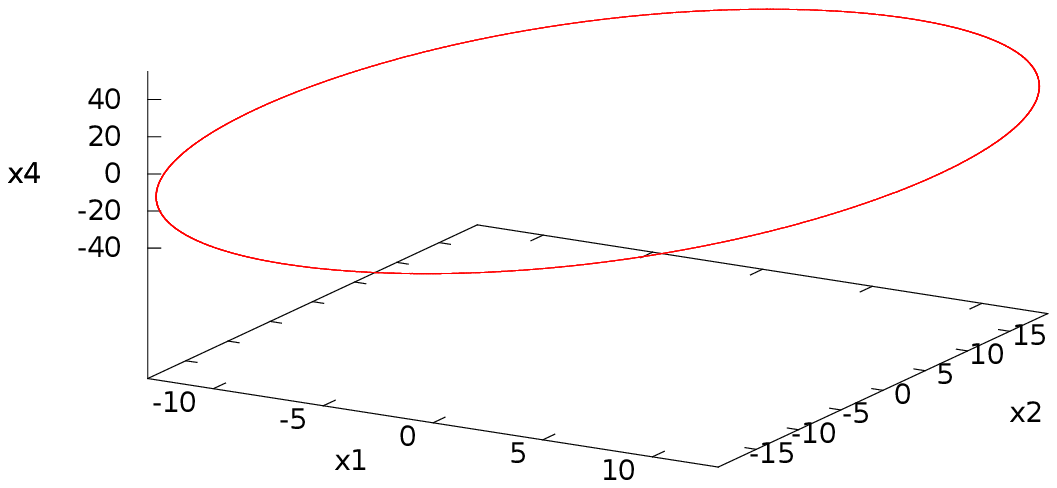}
   \caption{Projection of the cycle of the system \reff{dynsys} with the period $T_c \approx 0.3655$
   	        into the space ${x_1-x_2-x_4}$.}\label{fig:2}
\end{figure}

\begin{figure}[!]
   \centering\includegraphics[width=\textwidth,height=9cm]{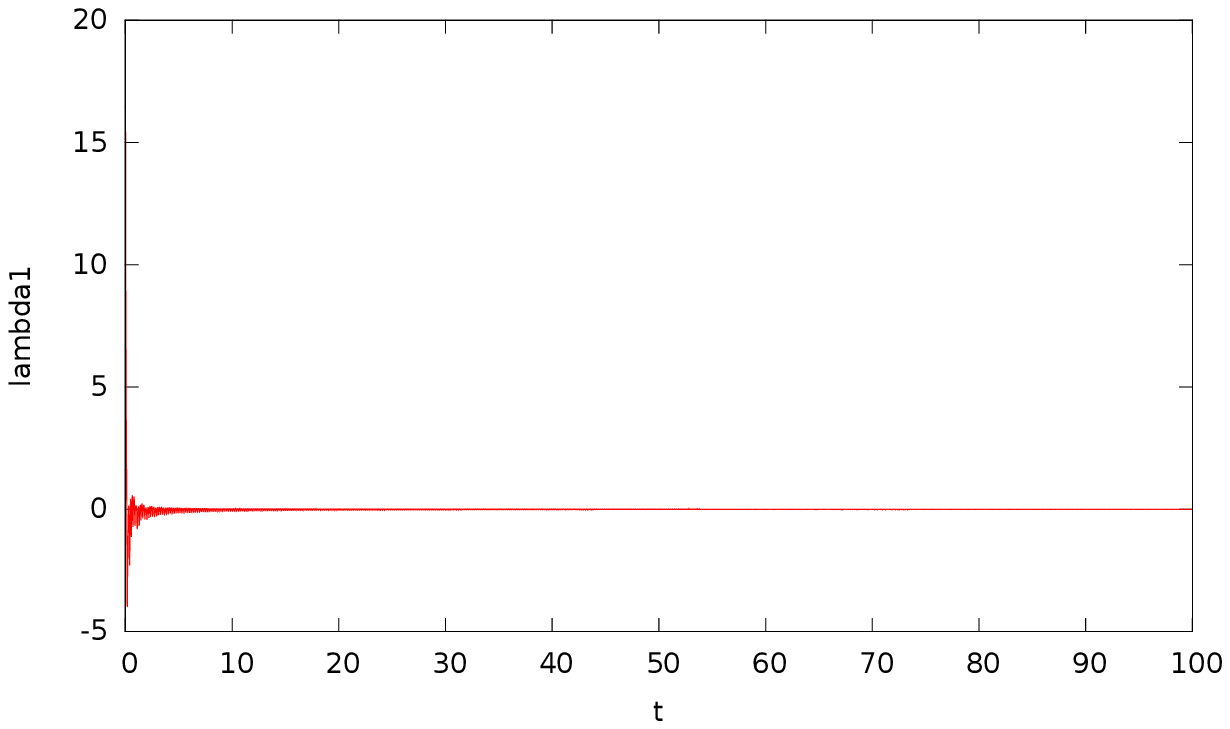}
   \caption{The graph of the change in time of the Lyapunov exponent $LE_1$.}\label{fig:3}
\end{figure}
\begin{figure}[!]
   \centering\includegraphics[width=\textwidth,height=9cm]{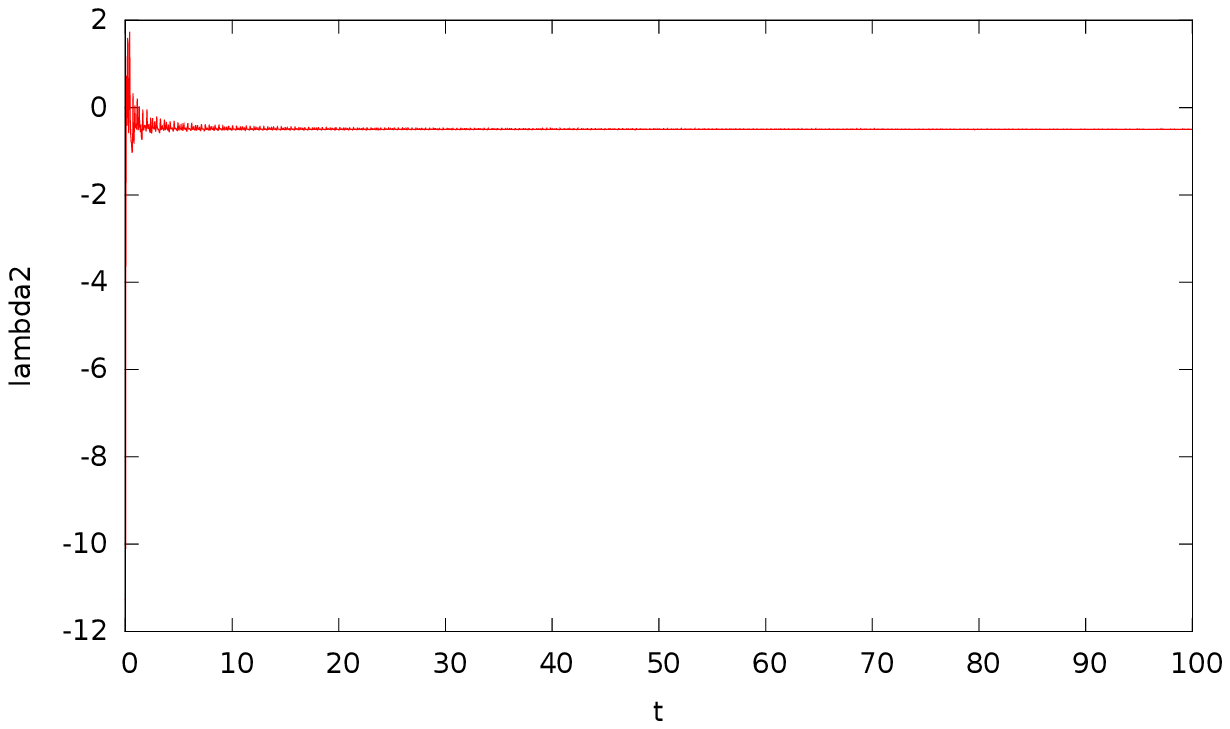}
   \caption{The graph of the change in time of the Lyapunov exponent $LE_2$.}\label{fig:4}
\end{figure}
\begin{figure}[!]
   \centering\includegraphics[width=\textwidth,height=9cm]{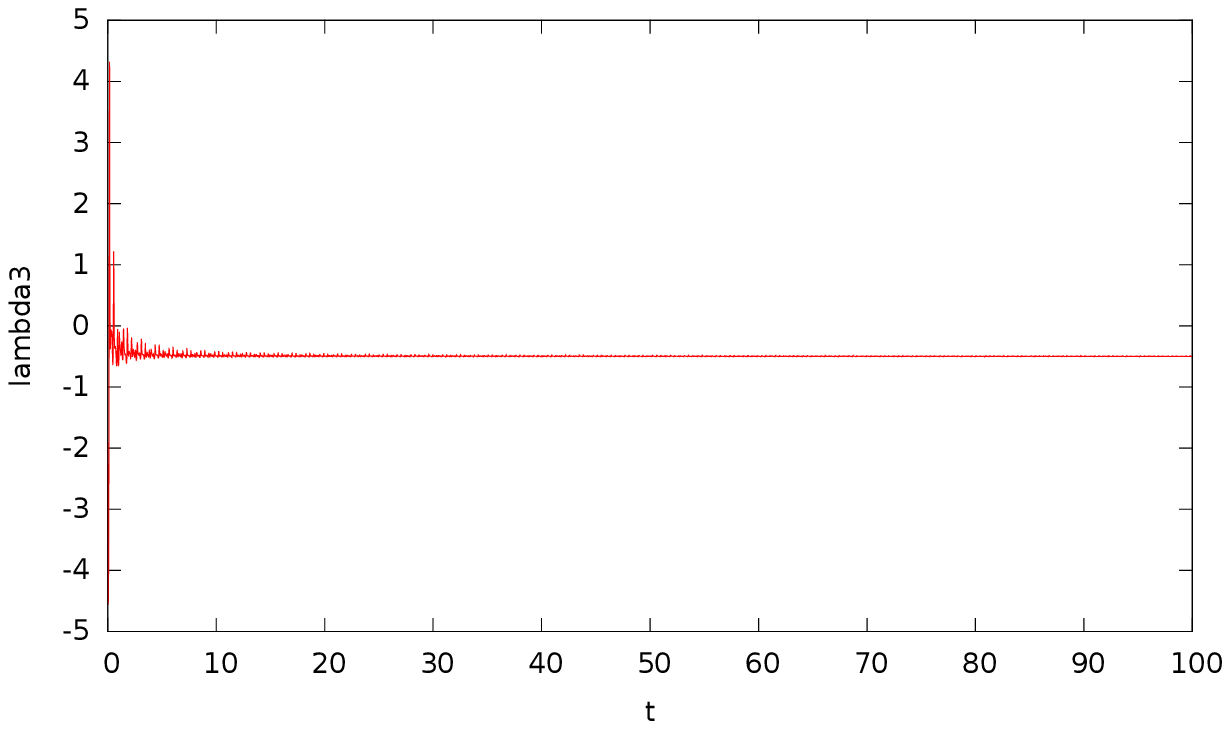}
   \caption{The graph of the change in time of the Lyapunov exponent $LE_3$.}\label{fig:5}
\end{figure}
\begin{figure}[!]
   \centering\includegraphics[width=\textwidth,height=9cm]{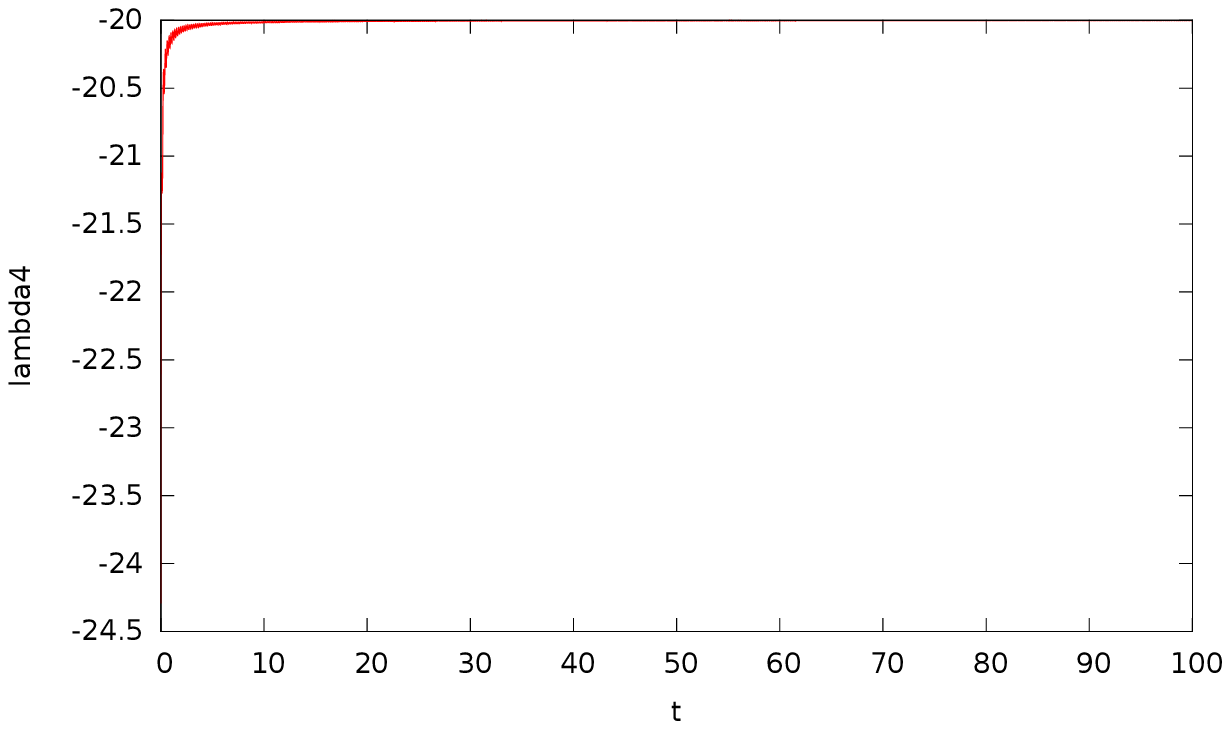}
   \caption{The graph of the change in time of the Lyapunov exponent $LE_4$.}\label{fig:6}
\end{figure}

From this iterative procedure, we see that the limiting trajectory in this case is a cycle with a period
$T_c \approx 0.3655$. Its projections are shown in Figs. \ref{fig:1} and \ref{fig:2}.

To select the value $\Delta t_P$, this step must be reduced until the rapprochement distance $d_{k^*}$ changes significantly.

This approach was also applied for other initial conditions. As a result, we get the same cycle.

Next, we investigate the stability of the resulting cycle by calculating the Lyapunov exponents.

\section{Calculation of the Lyapunov exponents}

In \cite{Pchel2020}, a modification of the algorithm for calculating the Lyapunov exponents
for a model describing the growth of a cancerous tumor was considered. Next, we will consider a generalization
of this algorithm for dynamical systems with a quadratic right-hand side of the $n$-th order.
Note that an other algorithm for calculating the Lyapunov exponents based on calculating
the fundamental matrix of the dynamical system and its QR decomposition is proposed in \cite{NVKuznetsov}.
The procedure described below can also be applied to it.

First, we define the form of the linearized system of equations. Let us denote the perturbations by
$$
  x_{n+1}(t),\:x_{n+2}(t),\:\ldots,\:x_{2n}(t).
$$

We define the form of the following expression using the rule for differentiating the sums and products:
$$
  \begin{array}{c}
     \dfrac{\partial\langle Q_pX,X \rangle}{\partial x_j}\cdot x_{n+j}=
     \dfrac{\partial}{\partial x_j}\left(\disp\sum_{i=1}^n \sum_{k=1}^n
     q_{i,k}^{(p)}x_i x_k\right)\cdot x_{n+j}=\nl
     =\left(\disp\sum_{k=1}^n q_{j,k}^{(p)}x_k+
     \sum_{i=1}^n q_{i,j}^{(p)}x_i\right)\cdot x_{n+j}=
     \disp\sum_{i=1}^n \left(q_{i,j}^{(p)}+q_{j,i}^{(p)}\right)x_i x_{n+j},
  \end{array}
$$
where $p=\overline{1;n}$, $j=\overline{1;n}$. Let us find the product
\equation
  \begin{array}{c}
  \left[\dfrac{\partial\langle Q_pX,X \rangle}{\partial x_1}\:\:
        \dfrac{\partial\langle Q_pX,X \rangle}{\partial x_2}\:\ldots\:
        \dfrac{\partial\langle Q_pX,X \rangle}{\partial x_n}\right]\cdot\left[
            \begin{array}{c}
               x_{n+1}\\
               x_{n+2}\\
               \vdots\\
               x_{2n}
            \end{array}
         \right]=\nl=
  \disp\sum_{j=1}^n\sum_{i=1}^n \left(q_{i,j}^{(p)}+q_{j,i}^{(p)}\right)x_i x_{n+j}=
  \disp\sum_{i=1}^n\sum_{j=n+1}^{2n}\left(q_{i,j-n}^{(p)}+q_{j-n,i}^{(p)}\right)x_i x_j.
  \end{array}\label{qu_form}
\endequation

Next, we extend the system \reff{dynsys} by adding to it a linearized system, which also
has a quadratic right-hand side. To do this, we introduce the extended matrix $\tilde A$:
$$
  \tilde A=\big[\tilde a_{i,j}\big]_{2n\times 2n},\:\:
  \tilde a_{i,j}=\left\{\begin{array}{l}
                           a_{i,j},\:\:\mbox{if}\:\:i=\overline{1;n}\:\:\&\:\:
                                                    j=\overline{1;n},\\
                           a_{i-n,j-n},\:\:\mbox{if}\:\:i=\overline{n+1;2n}\:\:\&\:\:
                                                    j=\overline{n+1;2n},\\
                           0\:\:\mbox{otherwise}.
                        \end{array}
                     \right.
$$
The vector function $\tilde X (t)$ is equal to
$$
  \tilde X(t)=\left[x_1(t)\:\:x_2(t)\:\ldots\:x_n(t)\:\:x_{n+1}(t)\:\ldots\:x_{2n}(t)\right]\tr.
$$

Based on the formula \reff{qu_form}, we have a quadratic form on the right-hand side of the
linearized system for each equation. Then, to reduce the system to general form, we introduce
the following extended matrix:
$$
  \begin{array}{l}
     \tilde Q_l=\left[\tilde q_{i,j}^{(l)}\right]_{2n\times 2n},\nl
     \mbox{for\:\:}l=\overline{1;n}:\:\:
     \tilde q_{i,j}^{(l)}=\left\{\begin{array}{l}
                                    q_{i,j}^{(l)},\:\:\mbox{if}\:\:i=\overline{1;n}\:\:\&\:\:
                                                                   j=\overline{1;n},\\
                                    0\:\:\mbox{otherwise},
                                 \end{array}
                          \right.\nl
     \mbox{for\:\:}l=\overline{n+1;2n}:\:\:
     \tilde q_{i,j}^{(l)}=\left\{\begin{array}{l}
                                    q_{i,j-n}^{(l-n)}+q_{j-n,i}^{(l-n)},\:\:\mbox{if}\:\:
                                                                   i=\overline{1;n}\:\:\&\:\:
                                                                   j=\overline{n+1;2n},\\
                                    0\:\:\mbox{otherwise}.
                                 \end{array}
                          \right.
  \end{array}
$$

We get an extended system similar to the system \reff{polsys}:
\equation
   \dot{\tilde X}=\tilde A \tilde X+\tilde\Phi(\tilde X),\label{expolsys}
\endequation
where
$$
\begin{array}{c}
   \tilde\Phi=\left[\tilde\varphi_1(\tilde X)\:\ldots\:\tilde\varphi_{2n}(\tilde X)\right]\tr,\nl
   \tilde\varphi_l(\tilde X)=\langle \tilde Q_l \tilde X,\tilde X \rangle,\:\:
   l=\overline{1;2n}.
\end{array}
$$

To obtain the domain of convergence of the power series that is a solution to the system \reff{expolsys},
the norms of the matrices are calculated
$$
  \|\tilde A\|,\:
  \|\tilde Q_1\|,\:
  \|\tilde Q_2\|,\:\ldots\:
  \|\tilde Q_{2n}\|
$$
and the number
$$
  \tilde\mu=2n\max_{p=\overline{1;2n}}\|\tilde Q_p\|.
$$

Let
\equation
   \tilde X(t)=\sum_{i=0}^\infty \tilde\lam_i t^i,\label{expowser}
\endequation
where $\tilde\lam_0=\tilde X(0)$ is the initial condition for the system \reff{expolsys},
$$
  \tilde\Psi_i=\left[\tilde\psi_{i,1}\:\ldots\:\tilde\psi_{i,2n}\right]\tr
$$
and
$$
  \tilde\psi_{i,p}=\prodser{i}{p}{\tilde},\:\:p=\overline{1;2n}.
$$

By analogy with the system \reff{polsys}, two auxiliary numbers are calculated as
functions of $\tilde \lam_0$:
$$
\begin{array}{c}
   \tilde h_1(\tilde\lam_0)=\|\tilde\lam_0\|,\nl
   \tilde h_2(\tilde\lam_0)=\left\{\begin{array}{l}
                                      \tilde\mu \tilde h_1^2(\tilde\lam_0)+(\|\tilde A\|+2\tilde\mu)
                                      \tilde h_1(\tilde\lam_0),
                                                 \:\:\mbox{if}\:\:\tilde h_1>1,\nl
                                      \|\tilde A\|+\tilde\mu\:\:\mbox{otherwise}.
                                   \end{array}
                            \right.
\end{array}
$$
Then the power series \reff{expowser} converges for $t\in[-\tilde\tau,\tilde\tau]$, where
$$
  \tilde\tau=\tilde\tau(\tilde\lam_0)=\dfrac{1}{\tilde h_2(\tilde\lam_0)+\delta},
$$
$\delta$ is any positive number.

Let us consider a modification of Benettin's algorithm for calculating the Lyapunov exponents using the
Gram--Schmidt orthogonalization \cite[pp. 163-165]{Kuznetsov}:

\begin{enumerate}
   \item \textbf{Divide} a given time interval $[0,T]$ (usually the value $T$ is large) into short intervals
         by the length
   $$
     \tau_M=\dfrac{T}{M},
   $$
   where $M$ is the number of such intervals.
   \item \textbf{Let} $X^*$ is a point of the researched solution; e.g., \reff{init_cond3};
   \item $k:=0$;
   \item $Y^{(k)}:=X^*$;
   \item \textbf{Let} $Z^{(k)}_{(p)}$ is the column vectors of initial perturbations from $n$-components
         (the real numbers), $p=\overline{1;n}$ ($p$ is a number of the perturbation vector);
   \item \textbf{Assign} to each component of the vector $Z^{(k)}_{(p)}$ random number in the range $[0,1]$;
   \item $LE_p:=0$, $p=\overline{1;n}$ (the initial values of the sums at calculating the Lyapunov exponents);
   \item \textbf{Orthogonalize} and \textbf{normalize} to unity the system of vectors $Z^{(k)}_{(1)}$,
         $Z^{(k)}_{(2)}$, ..., $Z^{(k)}_{(n)}$;
   \item\label{stepbc}$k:=k+1$;
   \item \textbf{For} $p$ from $1$ to $n$ with step $1$:\\
         \textbf{\textit{Beginning of cycle}}\\
         \textbf{Let}
         $$
           \tilde X(0):=\left[\begin{array}{r}
                                 Y^{(k-1)}\\
                                 Z^{(k-1)}_{(p)}
                              \end{array}
                        \right].
         $$
         \textbf{Find} $\tilde X(\tau_M)$ for the system \reff{expolsys} on algorithm, described in Sec.
         \ref{algo_FGBFI}.\\
         \textbf{Assign} to vectors $Y^{(k)}$ and $Z^{(k)}_{(p)}$ the matching vector blocks
         $\tilde X(\tau_M)$, i.e.
         $$
           \left[\begin{array}{r}
                    Y^{(k)}\\
                    Z^{(k)}_{(p)}
                 \end{array}
           \right]:=\tilde X(\tau_M).
         $$
         Note that for all values of $p$ the vector $Y^{(k)} $ is the same.\\
         \textbf{\textit{End of cycle}}
   \item \textbf{For} $p$ from $1$ to $n$ with step $1$:\\
         \textbf{\textit{Beginning of cycle}}\\
         $\disp S:=Z^{(k)}_{(p)}-\sum_{i=1}^{p-1}
                   \left\langle Z^{(k)}_{(p)},Z^{(k)}_{(i)}\right\rangle\cdot Z^{(k)}_{(i)};$\\
         $LE_p:=LE_p+\log|S|;$\\
         $Z^{(k)}_{(p)}:=\dfrac{S}{|S|}.$\\
         \textbf{\textit{End of cycle}}
   \item \textbf{If} $k<M$, \textbf{then} \textbf{Goto} step \ref{stepbc};
   \item $LE_p:=\dfrac{LE_p}{M\cdot\tau_M}$, $p=\overline{1;n}$;
   \item \textbf{Print} $LE_p$, $p=\overline{1;n}$.
\end{enumerate}

For the point \reff{init_cond3}, the Lyapunov exponents were found using the described algorithm for
$T=100$, $M=20000$ and $\tau_M=0.005$. Next, we give them the approximate values:
$$
  LE_1 \approx 0,\:\:LE_2 \approx -0.498,\:\:LE_3 \approx -0.499,\:\:LE_4 \approx -20.002.
$$
As can be seen, their signature corresponds to a stable limit cycle in the system \reff{dynsys},
which confirms the presence of the cycle, found through the study of the Poincar\'e recurrences.
Note that in \cite{PchelDiff} an alternative numerical-analytical scheme for constructing
periodic solutions of systems with a quadratic right-hand side was proposed using an example of
the Lorenz system. Now this scheme is actively developing in \cite{Luo1,Luo2,Luo3}
to find periodic solutions of nonlinear systems of differential equations.

The graphs of changes in the Lyapunov exponents over time are shown in Figs. \ref{fig:3} -- \ref{fig:6}.

The source codes of all computing programs are located at \cite{PchelProg}.

\section{Conclusion}

In this article, the author investigates the limiting regimes corresponding to different initial points of
the phase space. As a result, we get the same cycle shown in Fig. \ref{fig:1} and \ref{fig:2}. Consequently,
there is no hyperchaos in system \reff{dynsys}.

The analysis of the Poincar\'e recurrences, described in the goals of research, showed that
the rapprochement with the initial point occurs approximately at the same time (i.e., we have a regularity).
This indicates a limit cycle in the system \reff{dynsys}.

Stabilization of the behavior of the values of the Lyapunov exponents in Figs. \ref{fig:3} -- \ref{fig:6}
indicates the correctness of finding their limiting values over a relatively short time interval $[0, 100]$.

Thus, the limiting trajectory investigated using the analysis of the Poincar\'e recurrences and the
high-precision calculation of the Lyapunov exponents (i.e., different methods) is the limit cycle
which contradicts the results in \cite{Dong2019} for the values of the parameters $a=7$, $b=50$, $c=3$, $d=10$, $e=5$, $f=5$ and $k=1.5$ of the system \reff{dynsys}. The obtained data is available at \cite{PchelProg}. Most likely, the authors of the article \cite{Dong2019} made a typo indicating such values of the parameters at which they found a hyperchaos.

A significant novelty of the article is the high-precision method for searching for the Poincar\'e recurrences
on attractors of dynamical systems with a quadratic right-hand side, also a modification of Benettin's algorithm
for finding the Lyapunov exponents in order to check the found limiting regimes.

The methods described in the article can be used to study not only periodic regimes, but also quasiperiodic
and chaotic regimes for dynamical systems with a quadratic right-hand side. According to the statistics of
Poincar\'e recurrences, we can determine the dimension of the attractor and compare it with the dimension obtained
when calculating the Lyapunov exponents (for example, the articles \cite{Anishchenko2,Anishchenko3}).

The advantage of the considered algorithm of calculating the Lyapunov exponents is combination of the
linearized system of equations and the original dynamical system in general form \reff{expolsys} for finding
the state and perturbation vectors together.

\section{Acknowledgments}

The reported study was funded by RFBR for the research project No. 20-01-00347.

\end{paracol}

\reftitle{References}

\end{document}